\input amstex
\documentstyle {amsppt}
\pageheight{50.5pc}
\pagewidth{32pc}

\topmatter
\title{}
{ Poisson Approximation of Processes with Locally Independent
Increments with Markov Switching }
\endtitle
\author
V. S. Koroliuk $^1$, N. Limnios$^2$ and I.V. Samoilenko$^1$
\endauthor
\address
$^1$Institute of Mathematics, Ukrainian National Academy of Science,
Kiev, Ukraine $^2$Laboratoire de Math\'ematiques Appliqu\'ees,
Universit\'e de Technologie de Compi\`egne, France
\endaddress
\email isamoil\@imath.kiev.ua
\endemail
\abstract In this paper, the weak convergence of additive
functionals of processes with locally independent increments and
with Markov switching in the scheme of Poisson approximation is
proved. For the relative compactness, a method proposed by R.
Liptser for semimartingales is used with a modification, where we
apply a solution of a singular perturbation problem instead of an
ergodic theorem.
\endabstract

\endtopmatter
\rightheadtext{ Poisson Approximation of PLII with Markov Switching}
\leftheadtext{ V. S. Koroliuk, N. Limnios and I.V. Samoilenko }

\document

\head
1. Introduction
\endhead
Poisson approximation is still an active area of research in several
theoretical and applied directions. Several recent works on this
topic can be found in the literature: we can find the classical
approach in \cite{1,2,3}, and the functional approach in
\cite{8,9,7,12}.

In particular in \cite{8,9} it has been studied the following
stochastic additive functional
$$
\xi(t)=\xi_0+\int_0^t\eta(ds;{x}(s)),\quad t\ge 0,\tag1$$ of a jump
Markov process with locally independent increments (PLII) (\cite{8,
page 14}) $\eta(t;\cdot), t\ge0$, (also known as a Piecewise
deterministic Markov process -- PDMP, \cite{5, Chapter 2}),
perturbed by the jump Markov process $x(t), t\ge 0$. The process (1)
is studied in a (functional) Poisson approximation scheme, within an
{\it ad hoc} time-scaling as we can see below (2).

In the Poisson approximation scheme, the jump values of the
stochastic system are split into two parts: a small jump taking
values with probabilities close to one and a big jump taken values
with probabilities tending to zero together with the series
parameter $\varepsilon\to 0$. So, in the Poisson approximation
principle the probabilities (or intensities) of jumps are normalized
by the series parameter $\varepsilon >0$. Hence the time-scaled
family of processes $\xi^\varepsilon(t), t\ge 0, \varepsilon>0$, has
to be considered.

However, the method used here to prove the weak convergence is quite
different from the method proposed by other authors
(\cite{6}-\cite{17}): the main point is to prove convergence of
predictable characteristics of semimartingales which are integral
functionals of some switching Markov processes. But the main
difficulty is that the predictable characteristics of semimartingale
themselves depend on the process we study. Thus, to prove the
convergence of the process we should prove convergence of
predictable characteristics that depend on the process. Ordinary
methods can't help in this situation separately.

We propose to study functionals of PLII \cite{8, page 14} using a
combination of two methods. The method proposed by R. Liptser in
\cite{11}, based on semimartingales theory, is combined with a
solution of singular perturbation problem instead of ergodic
theorem. So, the method includes two steps.

At the first step we prove the relative compactness of the
semimartingale representation of the family $\xi^\varepsilon$,
$\varepsilon>0$, by proving the following two facts as proposed in
Liptser \cite{11}:
$$\lim\limits_{c\to \infty}\sup\limits_{\varepsilon\leq\varepsilon_0}
\Bbb P\{\sup\limits_{t\leq T}|\xi^{\varepsilon}(t)|>c\}=0, \hskip
2mm \forall \varepsilon_0>0$$ that is known as compact containment
condition (CCC), and
$$\Bbb E|\xi^{\varepsilon}(t)-\xi^{\varepsilon}(s)|^2\le k
|t-s|,$$ for some positive constant $k$.

At the second step we prove convergence of predictable
characteristics of the semimartingales, which are integral
functionals of the form ($a(u,x)$ is a real-valued function):
$$\int_0^t a(\xi^\varepsilon(s), x^\varepsilon(s))ds,$$
by using singular perturbation technique as presented in \cite{8}.

Finally, we apply Theorem IX.3.27 from \cite{7} in order to prove
the weak convergence of semimartingale.

The paper is organized as follows. In Section 2 we present the
time-scaled additive functional (1), the PLII and the switching
Markov process. In the same section we present the main results of
Poisson approximation. In Section 3 we present the proof of the
theorem.

\head 2. Main results
\endhead

Let us consider the space $\Bbb R^d$ endowed with a norm $|\cdot|$
($d\ge 1$), and $(E,\Cal E)$, a {\it standard phase space}, (i.e.,
$E$ is a Polish space and $\Cal E$ its Borel $\sigma$-algebra). For
a vector $v\in \Bbb R^d$ and a matrix $c\in \Bbb R^{d\times d}$ ,
$v^*$ and $c^*$ denote their transpose respectively. Let $C_3(\Bbb
R^d)$ be a measure-determining class of real-valued bounded
functions $g$ such that $g(u)/|u|^2 \to 0$, as $|u|\to 0$ (see
\cite{7,8}).

The additive functional $\xi^{\varepsilon}(t), t\geq 0,
\varepsilon>0$ on $\Bbb R^d$ in the series scheme with small series
parameter $\varepsilon\to 0$, $(\varepsilon>0)$ is defined by the
stochastic additive functional \cite{8, Section 3.3.1}
$$\xi^\varepsilon(t)=\xi_0^\varepsilon
+\int_0^t\eta^{\varepsilon}(ds;{x}(s/\varepsilon)).\tag2$$

The family of the Markov jump processes with {\it locally
independent increments} $\eta^{\varepsilon}(t;x),$ $t\geq 0, x\in E$
on $\Bbb R^d$, is defined by the generators on the test-functions
$\varphi(u)\in C^1(\Bbb R^d)$ \cite{8, Section 3.3.1} (see also
\cite{9})
$$\widetilde{\Gamma}^{\varepsilon}(x)\varphi(u)=\varepsilon^{-1}\int_{\Bbb R^d}[\varphi(u+v)
-\varphi(u)]\Gamma^{\varepsilon}(u,dv;x),\quad  x\in E,\tag3$$ or,
equivalently
$$\widetilde{\Gamma}^{\varepsilon}(x)\varphi(u)=b_\varepsilon(u;x)\varphi'(u)
+\frac{1}{2}c_\varepsilon(u;x)\varphi''(u)+\varepsilon^{-1}
\int_{\Bbb R^d}[\varphi(u+v)-\varphi(u)-v\varphi'(u)$$
$$-\frac{v^2}{2}\varphi''(u)]\Gamma^{\varepsilon}(u,dv;x),$$ where
$b_\varepsilon(u;x)=\varepsilon^{-1}\int_{\Bbb
R^d}v\Gamma^{\varepsilon}(u,dv;x)$,
$c_\varepsilon(u;x)=\varepsilon^{-1}\int_{\Bbb
R^d}vv^*\Gamma^{\varepsilon}(u,dv;x)$, and
$\Gamma^{\varepsilon}(u,dv;x)$ is the intensity kernel.

The switching Markov process ${x}(t), t\ge 0$ on the standard phase
space $(E,\Cal E)$, is defined by the generator
$$\Bbb Q\varphi(x) = q(x)\int_E
P(x,dy)[\varphi(y)-\varphi(x)],\tag4$$ where $q(x), x\in E$, is the
intensity of jumps function of $x(t), t\ge 0$, and $P(x,dy)$ the
transition kernel of the embedded Markov chain $x_n,n \ge 0$,
defined by $x_n= x(\tau_n), n\ge 0$, where $0=\tau_0\le \tau_1\le
...\le \tau_n\le ...$ are the jump times of $x(t), t\ge 0$. We
suppose also that the processes $\eta^\varepsilon(t;x)$ and $x(t)$
are right continuous.

It is worth noticing that the coupled process $\xi^\varepsilon(t),
x(t/\varepsilon), t\ge 0$ is a Markov additive process (see, e.g.,
\cite{8, Section 2.5}).

The Poisson approximation of Markov additive process (2) is
considered under the following conditions.

{\bf C1:} The Markov process ${x}(t), t \geq 0$ is uniformly ergodic
with $\pi(B), B\in \Cal E$ as a stationary distribution.

{\bf C2:} {\it Poisson approximation}. The family of processes with
locally independent increments $\eta^{\varepsilon}(t;x), t\geq 0,
x\in E$ satisfies the Poisson approximation conditions \cite{8,
Section 7.2.3}:

{\bf PA1:} Approximation of the mean values:
$$b_{\varepsilon}(u;x) = \int_{\Bbb R^d} v\Gamma^{\varepsilon}(u,dv; x)
= \varepsilon[b(u;x) +\theta_b^{\varepsilon} (u;x)],$$ and
$$c_{\varepsilon}(u;x) = \int_{\Bbb R^d}
vv^*\Gamma^{\varepsilon}(u,dv; x) = \varepsilon[c(u;x) +
\theta_c^{\varepsilon} (u;x)].$$

{\bf PA2:} Poisson approximation condition for intensity kernel
$$\Gamma_g^{\varepsilon}(u;x) = \int_{\Bbb R^d} g(v)\Gamma^{\varepsilon}(u,dv; x)
= \varepsilon[\Gamma_g(u;x) + \theta^{\varepsilon}_g(u;x)]$$ for all
$g \in C_3(\Bbb R^d)$, and the kernel $\Gamma_g(u;x)$ is bounded for
each $g \in C_3(\Bbb R^d)$, that is,
$$|\Gamma_g(u;x)| \leq\Gamma_g \quad \hbox{(a constant depending on $g$)}.$$

The above negligible terms
$\theta_a^\varepsilon,\theta_b^\varepsilon, \theta_c^\varepsilon$
satisfy the condition $$\sup\limits_{x\in E}
|\theta_{\cdot}^{\varepsilon}(u;x)|\to 0,\quad  \varepsilon\to 0.$$

In addition the following conditions are used:

{\bf C3:} {\it
Uniform square-integrability}:
$$\lim\limits_{c\to\infty}\sup\limits_{x\in E} \int_{|v|>c} vv^*\Gamma(u,dv; x) = 0,$$
where the kernel $\Gamma(u,dv; x)$ is defined on the class $C_3(\Bbb
R^d)$ by the relation $$\Gamma_g(u;x) = \int_{\Bbb R^d}
g(v)\Gamma(u,dv; x),\quad g \in C_3(\Bbb R^d).$$

{\bf C4:} {\it Linear growth}: there exists a positive constant $L$
such that
$$|b(u;x)|\leq L(1+|u|),\quad\hbox{and}\quad |c(u;x)|\leq L(1+|u|^2),$$
and for any real-valued non-negative function $f(x), x\in E$, such
that $$\int_{\Bbb R^d\setminus \{0\}}(1+f(x))|x|^2dx<\infty,$$ we
have
$$|\Lambda(u,v;x)|\leq Lf(v)(1+|u|),$$
where $\Lambda(u,v;x)$ is the Radon-Nikodym derivative of
$\Gamma(u,B;x)$ with respect to Lebesgue measure $dv$ in $\Bbb R^d$,
that is,
$$\Gamma(u,dv;x)=\Lambda(u,v;x)dv.$$

The main result of our work is the following.

\proclaim {Theorem 1} Under conditions {\bf C1-C4} the weak
convergence
$$\xi^{\varepsilon}(t)\Rightarrow \xi^0(t),\quad \varepsilon \to 0$$ takes
place.

The limit process $\xi^0(t), t\geq0$ is defined by the generator
$$\overline{\Gamma} \varphi(u)=\widehat{b}(u)\varphi'(u)+\int_{\Bbb
R^d}[\varphi(u+v)-
\varphi(u)-v\varphi'(u)]\widehat{\Gamma}(u,dv),\tag5$$ where the
average deterministic drift is defined by
$$\widehat{b}(u)=\int_E\pi(dx)b(u;x),$$
and the average intensity kernel is defined by
$$\widehat{\Gamma}(u,dv)=\int_E\pi(dx)\Gamma(u,dv;x).$$
\endproclaim

\remark{Remark 1} The limit process $\xi^0(t), t\ge 0$, is a  PLII
(see, e.g., \cite{8, page 14}), (or a PDMP - see, e.g., \cite{5,
Chapter 2}). The generator (5) can be written also as follows
$$\overline{\Gamma}\varphi(u)=\widehat{b}_0(u)\varphi'(u)+\int_{\Bbb
R^d}[\varphi(u+v)- \varphi(u)]\widehat{\Gamma}(u,dv),$$ where
$\widehat{b}_0(u)=\widehat{b}(u) - \int_{\Bbb R^d} v
\widehat{\Gamma}(u,dv)$.
\endremark

In the following corollary of the above theorem we give an important
particular case, where the limit process is a compound Poisson
process.

\proclaim{Corollary 1} Under Poisson approximation conditions:

PA1': Approximation of mean values:
$$b_{\varepsilon}(u;x) = \int_{\Bbb R^d} v\Gamma^{\varepsilon}(dv; x)
= \varepsilon [b(x)+ \theta^\varepsilon_b(u;x)]$$ and
$$c_{\varepsilon}(u;x) = \int_{\Bbb R^d} vv^*\Gamma^{\varepsilon}(dv; x)
= \varepsilon[c(x) + \theta_c^\varepsilon (u;x)].$$

PA2': Approximation condition for intensity kernel:
$$\Gamma_g^{\varepsilon}(u;x) = \int_{\Bbb R^d} g(v)\Gamma^{\varepsilon}(u,dv; x)
= \varepsilon[\Gamma_g(x) + \theta^{\varepsilon}_g(u;x)]$$ and the
kernel $\Gamma_g(x)$ is bounded for each $g \in C_3(\Bbb R^d)$, that
is,
$$|\Gamma_g(x)| \leq\Gamma_g \quad \hbox{(a constant depending on $g$)}.$$

And the additional condition

PA3:
$$\int_{\Bbb R^d} v\Gamma (dv)=\int_E \pi(dx)b(x),\quad \Gamma(dv)
=\int_{E}\pi(dx)\Gamma(dv;x),$$ the limit process $\xi^0(t), t\ge 0$
is a compound Poisson process
$$\xi^0(t)=u+ \sum_{k=1}^{\nu(t)} \alpha_k,\quad t\ge 0,$$
defined by the generator
$$\widetilde{\Gamma} \varphi (u)=\int_{\Bbb R^d} [\varphi(u+v)-\varphi(u)]\Gamma(dv),$$
where
$$\Gamma(dv)= \int_E \pi(dx)\Gamma(dv;x),\quad  \Gamma_g(x)= \int_{\Bbb R^d} g(v)\Gamma(dv;x).$$
The sequence of random variables $\alpha_k$, $k=1,2,...$ is i.i.d.
with joint distribution function $\Bbb P(\alpha_k\in dv) =
\Gamma(dv)/\Lambda, \quad \Lambda=\Gamma(\Bbb R^d)$ (it is obvious
that $\Gamma(\Bbb R^d)=\int_{\Bbb R^d}\Gamma(dv)$). The
time-homogeneous Poisson process $\nu (t), t\ge 0$, is defined by
its intensity: $\Lambda >0$.
\endproclaim

\head 3. Proof of Theorem 1
\endhead

The proof of Theorem 1 is based on the semimartingale representation
of the additive functional process (2). According to Theorems 6.27
and 7.16 \cite{4} the predictable characteristics of the
semimartingale (2) have the following representations:

$\bullet B^{\varepsilon}(t)=\varepsilon^{-1}\int_0^t
b_\varepsilon(\xi^{\varepsilon}(s);{x}^{\varepsilon}_s)ds= \int_0^t
b(\xi^{\varepsilon}(s);{x}^{\varepsilon}_s)ds+\theta^{\varepsilon}_{b},$

$\bullet C^{\varepsilon}(t)=\varepsilon^{-1}\int_0^t
c_\varepsilon(\xi^{\varepsilon}(s);{x}^{\varepsilon}_s)ds= \int_0^t
c(\xi^{\varepsilon}(s);{x}^{\varepsilon}_s)ds+\theta^{\varepsilon}_{c},$

$\bullet
\Gamma^{\varepsilon}(t)=\varepsilon^{-1}\int_0^t\int_{\Bbb{R}^d}g(v)\Gamma^{\varepsilon}
(\xi^{\varepsilon}(s),dv;{x}^{\varepsilon}_s)ds=\int_0^t\int_{\Bbb{R}^d}g(v)\Gamma
(\xi^{\varepsilon}(s),dv;{x}^{\varepsilon}_s)ds+\theta^{\varepsilon}_{g},$

where ${x}_t^{\varepsilon}:={x}(t/\varepsilon), t\ge 0$, and
$\sup\limits_{x\in E} |\theta_{\cdot}^{\varepsilon}|\to 0,\quad
\varepsilon\to 0.$

The jump martingale part of the semimartingale (2) is represented as
follows
$$\mu^{\varepsilon}(t)=\int_0^t\int_{\Bbb R^d}v[\mu^{\varepsilon}(\xi^{\varepsilon}(s),ds,dv;{x}^{\varepsilon}_s)-
\Gamma^{\varepsilon}(\xi^{\varepsilon}(s),dv;{x}^{\varepsilon}_s)ds].$$
Here $\mu^{\varepsilon}(u,ds,dv;x), x\in E$ is the family of
counting measures with characteristics
$$\Bbb E\mu^{\varepsilon}(u,ds,dv;x)=\Gamma^{\varepsilon}(u,dv;x)ds.$$

We can see now that predictable characteristics depend on the
process $\xi^{\varepsilon}(s).$ Thus, to prove convergence of
$\xi^{\varepsilon}(s)$ we should prove convergence of predictable
characteristics dependent on $\xi^{\varepsilon}(s).$ To avoid this
difficulty, we combine two methods.

We split the proof of Theorem 1 in the following two steps.

{\rm STEP 1}. At this step we establish the relative compactness of
the family of processes $\xi^{\varepsilon}(t), t\geq 0,
\varepsilon>0$ by using the approach developed in \cite{11}. Let us
remind that the space of all probability measures defined on the
standard space $(E,{\Cal E})$ is also a Polish space; so the
relative compactness and tightness are equivalent.

First we need the following lemma.

\proclaim {Lemma 1} Under assumption {\bf C4} there exists a
constant $k_T>0$, independent of $\varepsilon$ and dependent on $T$,
such that
$${\Bbb E}\sup\limits_{t\leq T}|\xi^{\varepsilon}(t)|^2\leq k_T.$$
\endproclaim

{\it Proof}: (following \cite{11}). The semimartingale (2) has the
following representation
$$\xi^{\varepsilon}(t)=u+A_t^{\varepsilon}+M_t^{\varepsilon},\tag6$$
where $u= \xi^\varepsilon(0)$; $A_t^{\varepsilon}$ is the
predictable drift
$$A_t^{\varepsilon}=\int_0^t b(\xi^{\varepsilon}(s),{x}_s^{\varepsilon})ds+
\int_0^t\int_{\Bbb{R}^d\setminus\{0\}}v\Gamma(\xi^{\varepsilon}(s),dv;{x}_s^{\varepsilon})ds+\theta^{\varepsilon},$$
and $M_t^{\varepsilon}$ is the locally square integrable martingale
$$M_t^{\varepsilon}=\int_0^t c(\xi^{\varepsilon}(s);{x}_s^{\varepsilon})dw_s+\int_0^t\int_{\Bbb{R}^d\setminus\{0\}}v[\mu^{\varepsilon}(ds,dv;{x}_s^{\varepsilon})
-\Gamma^{\varepsilon}(\xi^{\varepsilon}(s),dv;{x}_s^{\varepsilon})ds]+\theta^{\varepsilon},$$
where $w_t,t\ge 0$ is a standard Wiener process.

For a process $y(t), t\ge 0$, let us define the process
$y_{t}^\dag=\sup\limits_{s\leq t}|y(s)|,$ then from (6) we have
$$((\xi^{\varepsilon}_t)^\dag)^2\le
3[u^2+((A^{\varepsilon}_t)^\dag)^2+((M^{\varepsilon}_t)^\dag)^2].\tag7$$

Condition {\bf C4} implies that
$$(A^{\varepsilon}_t)^\dag \leq L\int_0^t(1+(\xi^{\varepsilon}_s)^\dag)ds
+\int_0^t\int_{\Bbb{R}^d\setminus\{0\}}|v|f(x)(1+(\xi^{\varepsilon}_s)^\dag)ds
$$ $$\leq L(1+r_1)\int_0^t(1+(\xi^{\varepsilon}_s)^\dag)ds,\tag8$$
where $r_1=\int_{\Bbb R^d\setminus \{0\}}|x|^2f(x)dx.$

Now, by Doob's inequality (see, e.g., \cite{12, Theorem 1.9.2}),
$$\Bbb E((M_t^{\varepsilon})^\dag)^2\leq
4|\Bbb E\langle M^{\varepsilon}\rangle_t|,$$ and condition {\bf C4}
we obtain
$$|\langle M^{\varepsilon}\rangle_t|=\left|\int_0^t
c(\xi^{\varepsilon}(s);{x}_s^{\varepsilon})
c^*(\xi^{\varepsilon}(s);{x}_s^{\varepsilon})ds+ \int_0^t\int_{\Bbb
R^d\setminus
\{0\}}vv^*\Gamma^{\varepsilon}(\xi^{\varepsilon}(s),dv;{x}_s^{\varepsilon})ds+\theta^{\varepsilon}\right|$$
$$\leq
 2L(1+r_1)\int_0^t[1+((\xi^{\varepsilon}_s)^\dag)^2]ds.\tag9$$

Inequalities (7)-(9) and Cauchy-Bunyakovsky-Schwartz inequality,
$([\int_0^t\varphi(s)ds]^2$ \hskip 2mm $\leq
t\int_0^t\varphi^2(s)ds)$, imply
$$\Bbb E((\xi^{\varepsilon}_t)^\dag)^2\leq
k_1+k_2\int_0^t\Bbb E((\xi^{\varepsilon}_s)^\dag)^2ds,$$ where $k_1$
and $k_2$ are positive constants independent of $\varepsilon$.

By Gronwall inequality (see, e.g., \cite{6, page 498}), we obtain
$$\Bbb E((\xi^{\varepsilon}_t)^\dag)^2\leq k_1\exp(k_2 t).$$

Hence the lemma is proved.

\proclaim {Corollary 2}  Under assumption {\bf C4}, the following
CCC holds:
$$\lim\limits_{c\to \infty}\sup\limits_{\varepsilon\leq\varepsilon_0}
\Bbb P\{\sup\limits_{t\leq T}|\xi^{\varepsilon}(t)|>c\}=0, \hskip
2mm \forall \varepsilon_0>0.$$
\endproclaim
{\it Proof}: The proof of this corollary follows from Kolmogorov's
inequality.

\remark{Remark 2} Another way to prove CCC is proposed in \cite{8,
Theorem 8.10} and used by other authors \cite{6, 17}. They use the
function $\varphi(u)=\sqrt{1+u^2}$ and prove corollary for
$\varphi(\xi^{\varepsilon}_t)$ by applying the martingale
characterization of the Markov process.

That can be easily proved, due to specific properties of
$\varphi(u).$
\endremark

\proclaim{Lemma 2} Under assumption {\bf C4} there exists a constant
$k>0$, independent of $\varepsilon$ such that
$$\Bbb E|\xi^{\varepsilon}(t)-\xi^{\varepsilon}(s)|^2\leq k |t-s|.$$
\endproclaim

\noindent{\it Proof}: In the same manner with (7), we may write
$$|\xi^{\varepsilon}(t)-\xi^{\varepsilon}(s)|^2\leq 2|A_t^{\varepsilon}
-A_s^{\varepsilon}|^2+2|M_t^{\varepsilon}-M_s^{\varepsilon}|^2.$$ By
using Doob's inequality, we obtain
$$\Bbb E|\xi^{\varepsilon}(t)-\xi^{\varepsilon}(s)|^2\leq
2\Bbb E\{|A_t^{\varepsilon}-A_s^{\varepsilon}|^2+8|\langle
M^{\varepsilon}\rangle_t-\langle M^{\varepsilon}\rangle_s|\}.$$

Now (8), (9), and assumption {\bf C4} imply
$$|A_t^{\varepsilon}-A_s^{\varepsilon}|^2+8|\langle
M^{\varepsilon}\rangle_t-\langle M^{\varepsilon}\rangle_s|\leq
k_3[1+((\xi_T^{\varepsilon})^\dag)^2]|t-s|,$$ where $k_3$ is a
positive constant independent of $\varepsilon$.

From the last inequality and Lemma 1 the desired conclusion emerges.

Thus from Corollary 2 and Lemma 2 immediately follows compactness of
the family of processes $\xi^{\varepsilon}(t), t\geq 0,
\varepsilon>0$.

{\rm STEP 2.} The next step of proof concerns the convergence of the
predictable characteristics. To do that, we apply the results of
Sections 3.2-3.3 in \cite{8} and the following theorem. $C^2_0(\Bbb
R^d\times E)$ is the space of real-valued twice continuously
differentiable by the first argument functions, defined on $\Bbb
R^d\times E$ and vanishing at infinity, and $C(\Bbb R^d\times E)$ is
the space of real-valued continuous bounded functions defined on
$\Bbb R^d\times E$.

\proclaim{Theorem 2} (\cite{8, Theorem 6.3}) Let the following
conditions hold for a family of coupled Markov process
$\xi^{\varepsilon}(t), x^{\varepsilon}(t), t\geq0, \varepsilon>0$:

{\bf CD1:} There exists a family of test functions
$\varphi^{\varepsilon}(u, x)$ in $C^2_0(\Bbb R^d\times E)$, such
that
$$\lim\limits_{\varepsilon\to 0}\varphi^{\varepsilon}(u, x) =
\varphi(u),$$ uniformly by $u, x.$

{\bf CD2:} The following convergence holds for the generator $\Bbb
L^{\varepsilon}$ of a coupled Markov process $\xi^{\varepsilon}(t),
x^{\varepsilon}(t), t\geq0, \varepsilon>0$
$$\lim\limits_{\varepsilon\to
0}\Bbb L^{\varepsilon}\varphi^{\varepsilon}(u, x) = \Bbb
L\varphi(u),$$ uniformly by $u, x$. The family of functions $\Bbb
L^{\varepsilon}\varphi^{\varepsilon}, \varepsilon>0$ is uniformly
bounded, both $\Bbb L\varphi(u)$ and $\Bbb
L^{\varepsilon}\varphi^{\varepsilon}$ belong to $C(\Bbb R^d\times
E)$.

{\bf CD3:} The quadratic characteristics of the martingales that
characterize a coupled Markov process $\xi^{\varepsilon}(t),
x^{\varepsilon}(t), t\geq0, \varepsilon>0$ have the representation
$\left\langle \mu^{\varepsilon}\right\rangle_t = \int^t_0
\zeta^{\varepsilon}(s)ds,$ where the random functions
$\zeta^{\varepsilon}, \varepsilon> 0,$ satisfy the condition
$$\sup\limits_{0\leq s \leq T} \Bbb E|\zeta^{\varepsilon}(s)|\leq
c < +\infty.$$

{\bf CD4:} The convergence of the initial values holds and
$$\sup\limits_{\varepsilon>0}\Bbb E|\zeta^{\varepsilon}(0)|\leq C
< +\infty.$$

Then the weak convergence
$$\xi^{\varepsilon}(t)\Rightarrow \xi(t),\quad \varepsilon\to 0,$$
takes place.
\endproclaim

We consider the three component Markov process
$B^{\varepsilon}(t),\xi^{\varepsilon}(t),{x}_t^{\varepsilon}, t\ge
0$ which can be characterized by the martingale
$$\mu_t^{\varepsilon}=\varphi(B^{\varepsilon}(t),\xi^{\varepsilon}(t),{x}_t^{\varepsilon})-
\int_0^t\Bbb L^{\varepsilon}
\varphi(B^{\varepsilon}(s),\xi^{\varepsilon}(s),{x}_t^{\varepsilon})ds,$$
where its generator $\Bbb L^{\varepsilon}$ has the following
representation \cite{8}
$$
\Bbb L^{\varepsilon}=\varepsilon^{-1}\Bbb
Q+\widetilde{\Gamma}^{\varepsilon}+\Bbb B^{\varepsilon},\tag10$$
with $\widetilde{\Gamma}^{\varepsilon}$ given by (3), $\Bbb Q$ given
by (4), and $\Bbb
B^\varepsilon(u;x)\varphi(v)=b_\varepsilon(u;x)\varphi'(v).$

According to \cite{8, Theorem 7.3}, under the conditions {\bf C1-C3}
the limit generator for $\widetilde{\Gamma}^{\varepsilon},
\varepsilon\to 0$ has the form (5). However in order to prove the
convergence of predictable characteristics, it is sufficient to
study the action of the generator $\Bbb L^{\varepsilon}$ on test
functions of two variables $\varphi(v,x)$.

Thus, it has the representation
$$\Bbb L^{\varepsilon}\varphi(v,x)=[\varepsilon^{-1}\Bbb Q+\Bbb B]\varphi(v,x).\tag11$$
The solution of the singular perturbation problem at the test
functions
$\varphi^{\varepsilon}(v,x)=\varphi(v)+\varepsilon\varphi_1(v,x)$ in
the form $\Bbb L^{\varepsilon}\varphi^{\varepsilon} =\widehat{\Bbb
L}\varphi+\theta^{\varepsilon}\varphi$ can be found in the same
manner with Proposition 5.1 in \cite{8}. That is
$$\widehat{\Bbb L}=\widehat{\Bbb B},\tag12$$
where $\widehat{\Bbb B}\varphi(v)=\widehat{b}\varphi'(v).$

Similar results can be proved for two other predictable
characteristics.

Now Theorem 2 may be applied.

We see from (10) and (12) that the solution of singular perturbation
problem for $\Bbb L^{\varepsilon}\varphi^{\varepsilon}(u,v;x)$
satisfies the conditions {\bf CD1, CD2}. Condition {\bf CD3} of this
theorem implies that the quadratic characteristics of the
martingale, corresponding to a coupled Markov process, is relatively
compact. The same result follows from the CCC (see Corollary 2 and
Lemma 2) by \cite{7}. Thus, the condition {\bf CD3} follows from the
Corollary 2 and Lemma 2. As soon as $B^{\varepsilon}(0)=B^0(0),
\xi^{\varepsilon}(0)=\xi^0(0)$ we see that the condition {\bf CD4}
is also satisfied. Thus, all the conditions of above Theorem 2 are
satisfied, so the weak convergence $B^{\varepsilon}(t)\Rightarrow
B^0(t)$ takes place.

By the same reasoning we can show the convergence of the processes
$C^{\varepsilon}(t)$ and $\Gamma^{\varepsilon}(t).$

{\it The final step of the proof} is achieved now by using Theorem
IX.3.27 in \cite{7}. Indeed all the conditions of this theorem are
fulfilled.

As we have mentioned, the square integrability condition 3.24
follows from CCC (see \cite{7}). The strong dominating hypothesis is
true with the majoration functions presented in the Condition {\bf
C4}. Condition {\bf C4} implies the condition of big jumps for the
last predictable measure of Theorem IX.3.27 in \cite{7}. Conditions
iv and v of Theorem IX.3.27 \cite{7} are obviously fulfilled.

The weak convergence of predictable characteristics is proved by
solving the singularly perturbation problem for the generator (11).

The last condition (3.29) of Theorem IX.3.27 is also fulfilled due
to CCC proved in Corollary 2 and Lemma 2. Thus, the weak convergence
is true.

We can see now that the limit Markov process is characterized by the
following predictable characteristics
$$B^0(t)=\int_0^tb(\xi^0(s))ds,\quad C^0(t)=\int_0^tc(\xi^0(s))ds,\quad
\Gamma^0_g(t)=\int_0^t\Gamma_{g}(\xi^0(s))ds.$$ So, the limit Markov
process $\xi^0(t)$ can be expressed by the generator (5).

Theorem 1 is proved.

{\it Acknowledgements.} The authors thank University of Bielefeld
for hospitality and financial support by DFG project 436 UKR
113/80/04-07. We also thank the anonymous referee for useful
comments and remarks.\par

\enddocument